\documentclass[review]{article}
\usepackage{mathtools}

\usepackage{amsmath,amssymb}
\usepackage{graphicx,epstopdf}
\usepackage{hyperref}
\usepackage{amsthm}
\urldef{\mailsa}\path{anandkumarv11@gmail.com}
\urldef{\mailsb}\path{r_radhaiyer@cb.amrita.edu}

\usepackage{graphicx}
\usepackage{lineno,hyperref}
\usepackage{blindtext,graphicx}
\usepackage[absolute]{textpos}
\usepackage{fancyhdr}
\usepackage{multirow}
\usepackage[a4paper,total={6in,8in}]{geometry}

\begin{document}

\title{\textbf{Harmonic Index of Total Graphs of Some Graphs}}
\author{\textbf{V. Anandkumar$^{1*}$,\ Radha Rajamani Iyer$^{2}$ }}
 \date{\small{Department of Mathematics,\\EASA College of Engineering and Technology,\\ Coimbatore -641 112,INDIA.\\ Department of Mathematics,\\ Amrita Vishwa Vidyapeetham,\\ Coimbatore -641 112,INDIA.}}
 \maketitle

\begin{abstract}
Let $G=G(V,E)$ be a graph with $\mid V \mid$ vertices and $\mid E \mid$ edges and total graph, $T(G)$ is obtained from $G$. In This paper we have study the Harmonic index of total graph for standards graphs, bipartite graph of particular type, regular graph, Grid graph $G_{m,n}$ and Complete Binary tree.\\\\
\textbf{2010 Mathematics Subject Classification:} Primary 92E10, Secondary 94C15\\
\textbf{Keywords:}Total graph, Harmonic Index, Grid graph, Binary tree.
\end{abstract}
\textbf{1.Introduction}\\

Let $G(V,E)$ be a simple connected undirected graph with $\mid V \mid=n$ vertices, $\mid E\mid=m$ edges, $uv\in E(G)$ if $u$ and $v$ are adjacent in $G$. The $Randi\acute{c}$ index of graph $G$ is denoted by $R(G)$ and it is defined as
$R(G)=\sum\limits_{uv\in E(G)}(d_{u} d_{v})^{\frac{-1}{2}}[11]$, where $d_{u}$ is the degree of the vertex $u$ in $G$. Few related variants of the $Randi\acute{c}$  index has been introduced in last few years.Sum-connectivity index is one among them introduced by zhou and Trinajstic in[14,15] and it is defined as $X(G)=\sum\limits_{uv\in E(G)}(d_{u}+ d_{v})^{\frac{-1}{2}}[11]$. $X_{\alpha}(G)=\sum\limits_{uv\in E(G)}(d(u)+d(v))^{\alpha}$, where $\alpha $ is the real number. It has been found that $X_{\alpha}(G)$ and $R(G)$ correlate well between themselves and with $\Pi$-electron energy of benzenoid hydrocarbons. In this paper we study the another related variant of the  index namely the $Randi\acute{c}$ Harmonic index of a graph $G$ and is denoted by $H(G)$ and defined as
$H(G)=\sum\limits_{uv\in E(G)}\frac{2}{d(u)+d(v)}.$\\

This Index was first appeared in[5]. Zhong found the minimum and maximum values of the harmonic index for simple connected graph and trees, and characterized the corresponding extremal graph.\\\\
\textbf{2.Harmonic Index of Total Graphs of Some Graphs}\\

For every graph $G$ there exist a graph called the \textit{line-graph} of $G$ and is denoted by $L(G)$, whose points correspond in a one-to-one manner with the lines of G in such a way that two points in the $L(G)$ are adjacent if and only if the corresponding lines are adjacent in $G$. Vertices and edges are elements of a graph $G$. For any graph $G$, \textit{total graph}, denoted by $T(G)$, is obtained from $G$ whose vertex set $V(T(G))=V(G)\cup E(G)$ and two vertices are adjacent if and only if the corresponding elements are adjacent or incident in $G$. Total graph is introduced by Mehdi Behzad and Gary Chartrand in the year 1966. They have presented the characterization of graphs hamiltonian total graphs and of graph whose total graphs are Eulerian.\\\\
\textbf{Observation}\\

If $v^{\prime}$ is the point in the total graph of $G$ corresponding to point $v$ of $G$, then the degree of $v^{\prime}$, denoted by deg $v^{\prime}$ , equals 2 $deg \hspace*{.1cm} v$.
Where $x=v_{1}v_{2}$ is a line of $G$ and $x^{\prime}$ is the corresponding point of $T(G)$, then deg $x^{\prime}=deg v_{1}+deg v_{2}$.\\

In this paper Exact values of Harmonic index of total graph of standard graphs such as Path, Cycle, Gird Graph, Bipartite graph of a particular type and complete binary tree are obtained.\\\\
\textbf{Proposition 2.1.}
For any  graph $P_{n}$ of order $n\geq4$ , $H[T(P_{n})]=n-\frac{269}{420}$.\\\\
\textbf{Proposition 2.2.}
For any star graph $S_{n}, H[T(S_{n})]=(n-1)[\frac{2}{n+2}+\frac{3(n+1)}{2(3n-1)}]$.\\\\
\textbf{Theorem 2.3.}
Let $G$ be a graph with $\mid V\mid$=n and $\mid E\mid =m$ then \\$H(T(G))=\frac{1}{2}H(G)+\sum\limits^{frac{1}{2}M_{1}(G)-m}_{\begin{array}{l}uv \in E(G)\\vw\in E(G)\end{array}}\frac{2}{d(u)+2d(v)+d(w)}+ 8\sum\limits^{m}_{uv\in E(G)}[\frac{[d(u)+d(v)]}{3(d(u)+d(v))^{2}+4d(u)d(v)}].$
{\it{proof:}}
Total graph of a graph $G$ has $|E(T(G))|=3m+\frac{1}{2}\sum\limits_{v_{i}}d(v_{i})^2-m$\\
$\hspace*{2.7cm}=2m+\frac{1}{2}\sum\limits_{v_{i}}d(v_{i})^2$\\
each edge belongs any one of the following category.\\

(i) edges between $u^{\prime}$ and $v^{\prime}$, where $u^{\prime}$ and $v^{\prime}$ are vertices corresponding to the vertices $u$ and $v$ in $G$.

(ii) edges between $X^{\prime}$ and $Y^{\prime}$, where $X^{\prime}$ and $Y^{\prime}$ are vertices corresponding to the edges $X$ and $Y$ in $G$.

(iii) edges between $X^{\prime}$ and $u^{\prime}$, where $X^{\prime}$ and $u^{\prime}$ are vertices corresponding to the edge $X$ and the vertex $u$ respectively, $X=uv\in E(G)$.

$\therefore H[T(G)]=\sum\limits^{2m+\frac{1}{2}\sum\limits^{n}_{v_{i}\in G}d(v_{i})^2}_{u^{\prime}v^{\prime}\in E(T(G))} \frac{2}{d(u^{\prime})+d(v^{\prime})}$\\
$\hspace*{2.3cm}=\sum\limits^{m}_{u^{\prime}v^{\prime}\in E(T(G))}\frac{2}{d(u^{\prime})+d(v^{\prime})}+\sum\limits^{\frac{1}{2}\sum\limits^{n}_{v_{i}\in G}d(v_{i}^2)-m}_{X^{\prime}Y^{\prime}\in E(T(G))}\frac{2}{d(X^{\prime})+d(Y^{\prime})}+\sum\limits^{2m}_{X^{\prime}u^{\prime}\in E(T(G))}\frac{2}{d(X^{\prime})+d(u^{\prime})}.$\\\\
$\hspace*{2.3cm}=A+B+C$.\\

where $A=\sum\limits^{m}_{u^{\prime}v^{\prime}\in E(T(G))}\frac{2}{d(u^{\prime})+d(v^{\prime})}=\sum\limits^m_{uv\in E(G)}\frac{2}{2(d(u)+d(v))}$\\\\
$\hspace*{1.9cm}=\frac{1}{2}H(G).\hfill(1)$\\

where $B=\sum\limits^{\frac{1}{2}\sum\limits^{n}_{v_{i}\in G}d(v_{i}^2)-m}_{X^{\prime}Y^{\prime}\in E(T(G))}\frac{2}{d(X^{\prime})+d(Y^{\prime})}=\sum\limits^{\frac{1}{2}M_{1}(G)-m}_{XY\in E(L(G))}\frac{2}{d(X)+d(Y)+4}.$\\\\
$\hspace*{2cm}=\sum\limits^{\frac{1}{2}M_{1}(G)-m}_{\begin{array}{l}uv \in E(G)\\vw\in E(G)\end{array}}\frac{2}{d(u)+2d(v)+d(w)}.\hfill(2)$

where $C=\sum\limits^{2m}_{X^{\prime}u^{\prime}\in E(T(G))}\frac{2}{d(X^{\prime})+d(u^{\prime})}.$\\\\
$\hspace*{2cm}=\sum\limits^{2m}_{X=uv\in E(G)}\frac{2}{2d(u)+d(X)+2}$\\\\
$\hspace*{2cm}=\sum\limits^{m}_{uv\in E(G)}[\frac{2}{2d(u)+d(u)+d(v)}+\frac{2}{d(u)+2d(v)+d(v)}]$\\\\
$\hspace*{2cm}=\sum\limits^{m}_{uv\in E(G)}[\frac{2}{3d(u)+d(v)}+\frac{2}{d(u)+3d(v)}]$\\\\
$\hspace*{2cm}=\sum\limits^{m}_{uv\in E(G)}[\frac{8d(u)+8d(v)}{3d(u)^{2}+9d(u)d(v)+d(u)d(v)+3d(v)^{2}}]$\\\\
$\hspace*{2cm}=\sum\limits^{m}_{uv\in E(G)}[\frac{8[d(u)+d(v)]}{3d(u)^{2}+10d(u)d(v)+3d(v)^{2}}]$\\\\
$\hspace*{2cm}=8\sum\limits^{m}_{uv\in E(G)}[\frac{[d(u)+d(v)]}{3(d(u)+d(v))^{2}+4d(u)d(v)}].\hfill(3)$\\\\
From (1),(2) and (3),\\
$H(T(G))=\frac{1}{2}H(G)+\sum\limits^{M_{1}(G)-m}_{\begin{array}{l}uv \in E(G)\\vw\in E(G)\end{array}}\frac{2}{d(u)+2d(v)+d(w)}+ 8\sum\limits^{m}_{uv\in E(G)}[\frac{[d(u)+d(v)]}{3(d(u)+d(v))^{2}+4d(u)d(v)}].$\\\\
\textbf{Theorem 2.4.}
If $G$ is regular graph then $H[T(G)]=\frac{n(k+2)}{4}.$\\
{\it{proof:}}
If $G$ is a regular then $T(G)$ is also regular and it has $\frac{n}{2}(2+k)$ vertices and degree of all the vertices are the same. Then $|E(T(G))|=\frac{nk}{2}(k+2)$.\\
Hence, H[T(G)]=$\sum\limits^{\frac{nk}{2}(k+2)}_{u^{\prime}v^{\prime}\in E(T(G))}\frac{2}{d(u^{\prime})+d(v^{\prime})}$\\\\
here degree of every vertex in $T(G)$ is 2k.\\
therefore,\\$\hspace*{1.8cm}=\frac{nk}{2}(k+2)\frac{2}{2k+2k}$\\\\
$H[T(G)]=\frac{n(k+2)}{4}.$\\\\
\textbf{Corollary 2.4.1.}
Let $C_{n}$ be a cycle with $n$ vertices then $H[T(C_{n})]=n$.\\
{\it{proof:}}
 Total graph $T(C_{n})$ is the 4-regular graph of the size $4n$.\\\\
 $\therefore H[T(C_{n})]=\sum\limits^{4n}_{u\prime v\prime\in E[T(C_{n})]}\frac{2}{d(u\prime)+d(v\prime)}$\\\\
               $\hspace*{1.9cm}=4n\frac{2}{4+4}$\\
              $\hspace*{.3cm}H[T(C_{n})]=n$.\\\\
\textbf{Corollary 2.4.2.}
 Let $K_{n}$ be a complete graph then $H[T(K_{n})]=\frac{n(n+1)}{4}$.\\
{\it{proof:}}
$T(K_{n}) is the $(2(n-1))-regular graph of the size $\frac{n(n-1)(n+1)}{2}$.\\\\
$\therefore H[T(K_{n})]=\sum\limits^{n(n-1)(n+1)}_{u\prime v\prime\in E[T(K_{n})]}\frac{2}{d(u\prime)+d(v\prime)}$\\\\
           $\hspace*{1.6cm}=\frac{n(n-1)(n+1)}{2}(\frac{2}{2(n-1)+2(n-1)})$\\\\
$H[T(K_{n})]=\frac{n(n+1)}{4}.$\\\\
\textbf{Theorem 2.5.}
Let graph $G$ be a bipartite graph with bipartition $X$ and $Y$ for all $i$.
$\deg v_{i}=\begin{cases}k_{1}, & \text{if} \hspace*{.2cm} v_{i}\in X\\ k_{2}, & \text {if}\hspace*{.2cm}  v_{i}\in Y\end{cases},\\ then   H[T(G)]=\frac{mk_{1}^2+nk_{2}^2}{2(k_{1}+k_{2})}+\frac{4(mk_{1}+nk_{2})(k_{1}+k_{2})}{(3k_{1}+k_{2})(k_{1}+3k_{2})}$.\\
{\it{proof:}}
Let $G$ be a bipartite graph with partition $\mid X\mid =m$ and $\mid Y\mid =n$.\\
number of edges in $G$ is $\frac{mk_{1}+nk_{2}}{2}$\\
In the Total graph of $G$, $|V(T(G))|=m+n+\frac{mk_{1}+nk_{2}}{2}$\\\\
$|E(T(G))|=\frac{mk_{1}+nk_{2}}{2}+\frac{nk_{2}(k_{2}-1)}{2}+\frac{mk_{1}(k_{1}-1)}{2}+2(\frac{mk_{1}+nk_{2}}{2})$\\\\
$|E(T(G))|=\frac{1}{2}[nk_{2}(k_{2}+2)+mk_{1}(k_{1}+2)]$\\\\
$\hspace*{0.2cm}H[T(G)]=\sum\limits^{\frac{1}{2}[nk_{2}(k_{2}+2)+mk_{1}(k_{1}+2)]}_{u^{\prime}v^{\prime}\in E[T(G)]} \frac{2}{d(u^{\prime})+d(v^{\prime})}$\\\\
$\hspace*{1.7cm}=\frac{mk_{1}+nk_{2}}{2}(\frac{2}{2(k_{1}+k_{2})})+\frac{mk_{1}(k_{1}-1)}{2}(\frac{2}{2(k_{1}+k_{2})})+\frac{nk_{2}(k_{2}-1)}{2}(\frac{2}{2(k_{1}+k_{2})})+$\\\\
$\hspace*{1.8cm}\frac{mk_{1}+nk_{2}}{2}(\frac{2}{3k_{1}+k_{2}})+\frac{mk_{1}+nk_{2}}{2}(\frac{2}{k_{1}+3k_{2}})$\\\\
$\hspace*{1.7cm}=\frac{mk_{1}^2+nk_{2}^2}{2(k_{1}+k_{2})}+(mk_{1}+nk_{2})(\frac{1}{3k_{1}+k_{2}}+\frac{1}{k_{1}+3k_{2}})$\\\\
$\hspace*{1.7cm}=\frac{mk_{1}^2+nk_{2}^2}{2(k_{1}+k_{2})}+(mk_{1}+nk_{2})(\frac{4k_{1}+4k_{2}}{(3k_{1}+k_{2})(k_{1}+3k_{2})})$\\\\
$\therefore H[T(G)]=\frac{mk_{1}^2+nk_{2}^2}{2(k_{1}+k_{2})}+\frac{4(mk_{1}+nk_{2})(k_{1}+k_{2})}{(3k_{1}+k_{2})(k_{1}+3k_{2})}$.\\\\
\textbf{Corollary 2.5.1.}
let $K_{m,n}$ be a complete bipartite graph then\\\\ $H[T(K_{m,n})]=\frac{5mn}{6(m+n)^2+8mn}.$\\\\
{\it{proof:}}
$H[T(K_{m,n})]=mn\frac{2}{2n+2m}+\frac{mn(m+n-2)}{2}\frac{2}{(m+n+m+n)}+mn(\frac{2}{2n+m+n})$\\\\
$\hspace*{2.7cm}+mn\frac{2}{2m+m+n}$\\\\
$\hspace*{2.3cm}=\frac{mn}{m+n}+\frac{mn(m+n-2)}{2(m+n)}+\frac{2mn}{2n+m+n}+\frac{2mn}{2m+m+n}$\\\\
$\hspace*{2.3cm}=\frac{2mn+mn(m+n)-2mn}{2(m+n)}+\frac{2mn}{3n+m}+\frac{2mn}{3m+n}$\\\\
$\hspace*{2.3cm}=\frac{mn}{2}+\frac{2mn}{3n+m}+\frac{2mn}{3m+n}$\\\\
 $H[T(K_{m,n})]=\frac{5mn}{6(m+n)^{2}+8mn}.$\\\\
\textbf{Theorem 2.6.}
Let $G$ be a gird graph then\\ $H[T(G_{m,n})]=\frac{M_{1}(G)+mn}{2}+(m+n)\frac{99}{260}-\frac{58237}{12870}$ where $m,n\geq5$.\\
{\it{proof:}}
The total graph of gird graph has\\ $\hspace*{2cm}\mid V\mid=3mn-m-n$ and \\
$\hspace*{2cm} \mid E\mid =\frac{1}{2}M_{1}(G)+4mn-2(m+n)$.\\
The degree pairs of the number of edges $8,8,4,16,8,4(2m+2n-13),4(3m+3n-16),2(m+n-4),4,8(m+n-5),\frac{1}{2}M_{1}(G)+4mn-26(m+n)+84$ in the $T(G_{m,n})$ are (4,5), (4,6), (5,5), (5,6), (5,7), (6,6), (6,7), (6,8), (7,7), (7,8), (8,8) respectively.\\\\
$\therefore H[T(G_{m,n})]=8(\frac{2}{4+5})+8(\frac{2}{4+6})+4(\frac{2}{5+5})+16(\frac{2}{5+6})+8(\frac{2}{5+7})+4(2m+2n-13)(\frac{2}{6+6})+2(3m+3n-16)(\frac{2}{6+7})+2(m+n-4)(\frac{2}{6+8})+4(\frac{2}{7+7})+8(m+n-5)(\frac{2}{7+8})+(\frac{1}{2}M_{1}(G)+4mn-26(m+n))+84(\frac{2}{8+8})$\\

$\hspace*{2.2cm}=\frac{1}{16}M_{1}(G)+\frac{mn}{2}+(\frac{16}{9}+\frac{12}{5}+\frac{32}{11}-\frac{38}{3}-\frac{64}{13}-\frac{4}{7}+\frac{21}{2})+(m+n)(\frac{4}{3}+\frac{16}{13}+\frac{16}{15}-\frac{52}{16})$
$\hspace*{2.2cm}=\frac{1}{16}M_{1}(G)+\frac{mn}{2}+\frac{1188}{3120}(m+n)-(\frac{287539}{90090})$\\
$H[T(G_{m,n})]=\frac{1}{16}M_{1}(G)+\frac{mn}{2}+\frac{99}{260}(m+n)-(\frac{41077}{12870}).$\\

\textbf{Theorem 2.7.}
For any complete binary tree T of level $l\geq3,$ $H(T(T))=2^{l}(\frac{223}{120})-\frac{2274}{1485}.$\\
{\it{proof:}}
Let T be a binary tree with level l. number of vertices and number of edges in the total graph T(T) are $2^{l+2}-3$ and $9.2^{l}-11$ respectively.\\
ie., $|V[T(T)]|=2^{l+2}-3$ and $|E[T(T)]|=9.2^{l}-11$\\
 the number of edges, $2^{l}, 2^{l}, 2^{l-1}, 2, (2^{l+1}+2), 1, 6, (9.2^{l-1}-22$) have the degree pairs(ie., end degree vertices) in the total graph are (2,4), (2,6), (4,4), (4,5), (4,6), (5,5), (5,6) and (6,6) respectively. \\
$H[T(T)]=\sum\limits^{9.2^{l}-11}_{u^{\prime}v^{\prime}\in E(T(G))}\frac{2}{d(u^{\prime})+d(v^{\prime})}\\\\
\hspace*{2.1cm}=2^{l}(\frac{2}{2+4})+2^{l}(\frac{2}{2+6})+2^{l-1}(\frac{2}{4+4})+2(\frac{2}{4+5})+(2^{l+1}+2)(\frac{2}{4+6})+\\\\\hspace*{2.3cm}(\frac{2}{5+5})+6(\frac{2}{5+6})
+(9.2^{l-1}-22)(\frac{2}{6+6})\\\\
\hspace*{2.1cm}=2^{l}(\frac{2}{6})+2^{l}(\frac{2}{8})+2^{l-1}(\frac{2}{8})+2(\frac{2}{9})+2^{l+1}(\frac{2}{10})+2(\frac{2}{10})+\frac{2}{10}+6(\frac{2}{11})\\\\\hspace*{2.3cm}+9.2^{l-1}(\frac{2}{12})-22(\frac{2}{12})\\\\
\hspace*{2.1cm}=2^{l}(\frac{2}{6}+\frac{2}{8}+\frac{2}{16}+\frac{4}{10}+\frac{18}{24})+4(\frac{1}{9}+\frac{1}{10}+\frac{1}{20}+\frac{3}{11}-\frac{11}{12})\\\\
\hspace*{2.1cm}=2^{l-1}(\frac{40+30+15+48+90}{240})-\frac{2274}{1485}\\\\
\hspace*{2.1cm}=2^{l+1}(\frac{223}{240})-\frac{2274}{1485}\\\\
\therefore  H(T(T))=2^{l}(\frac{223}{120})-\frac{2274}{1485}.$\\\\
\textbf{Acknowledgement.} The authors are thankful to the anonymous referees for their useful comments.\\\\
\textbf{References}
\begin{enumerate}
\item G.H. Fath-Tabar, A.Hamzeh, and S.Hossein-Zadeh, $GA_{2}$ index of some graph operation, Filomat, 24(1):21-28, 2010.
\item Xinli Xu, Relationships between Harmonic index and other topological indices, applied mathematical Sciences, 6(41):2013-2018, 2012.
\item Ismael G.Yero and Juan A.Rodriguez-Velazquez, On the Randic index of corono product graph, International Scholarly Research notices, 2011.
\item M.Randic, On characterization of molecular branching, J.Am.Chem.Soc.97(1975)6609-6615.
\item L.Pogliani, From molecular connectivity indices to semi empirical connectivity terms: recent trends in graph theoretical descriptors,
    Chem.Rev.100(2000)3827-3858.
\item I. Gutman, B.Furtula(Eds.),Recent Results in the Theory of Randic Index, Univ.Kragujevac,2008.
\item B.Zhou, N.Trinajstic, On general sum-connectivity index, J.Math.Chem.47(2010)210-218.
\item B.Lucic, N.Trinajstic, B.Zhou, Comparison between the sum-connectivity index and product-connectivity index for benzenoid
    hydrocarbons,Chem.Phys.Lett.475(2009)146-148.
\item Z.Du,B.Zhou, N.Trinajstic, On the general sum-connectivity index of trees, appl.Math.Lett.24(2011)402-405.
\item Lingping Zhong, The harmonic index for graphs, applied Mathematics Letters 25(2012)561-566.
\item Hanyuan Deng, S.Balachandran, S.K.Ayyaswamy, Y.B.VenkataKrishnan, On the harmonic index and the chromatic number of a graph, Discrete applied
    Mathematics 161(2013)2740-2744.
\item Z.Du, B.Zhou, On sum-connectivity index of bicyclic graphs, Bull.Math.Sci.Soc.(2)35(1)(2012)101-117.
\item Z.Du, B.Zhou, N.Trinajstic, minimum sum-connectivity indices of trees and unicyclic graphs of a given matching number, J.Math.Chem.47(2010)842-855.
\item R.Xing, B.Zhou, N.Trinajstic, sum-connectivity index of molecular trees,J.Math.Chem.48(2010)583-591.
\item L.Zhong, The harmonic index for graphs, appl.Math.Lett.25(2012)561-566.
\item B.Shwetha Shetty, V.Lokesha and P.S.Ranjini, On the harmonic index of graph operations,Transactions on combinatorics,4(4)(2015),5-14.
\end{enumerate}

\end{document}